\documentclass[letterpaper, 10 pt, conference]{ieeeconf}

\IEEEoverridecommandlockouts                              % This command is only needed if 
\usepackage{amsmath} % assumes amsmath package installed
\usepackage{amssymb}  % assumes amsmath package installed
\usepackage[short, c3, nocomma]{optidef}
\usepackage{multirow}
\usepackage{graphicx}
\usepackage{graphicx}
\usepackage{amsmath}
\usepackage{amssymb}
\usepackage{amsthm}
\usepackage{dsfont}
\usepackage{commath}
\usepackage{mathtools}
\usepackage{algorithm}
\usepackage[noend]{algorithmic}
\usepackage{tabularx}
\usepackage{accents}

\usepackage[noadjust]{cite}

\usepackage{scalerel}

\usepackage[american]{circuitikz}
\usepackage{tikz}

\newtheorem{lem}{Lemma}

\theoremstyle{definition}
\newtheorem{defn}{Definition}

\newtheorem{rem}{Remark}
\newtheorem*{remark}{Remark}
\newtheorem{problem}{Problem}

\theoremstyle{definition}

\newcommand{\X}{\mathbf{\mathcal{X}}}

\newcommand{\editcolor}{black}
\makeatletter
\newcommand\fs@betterruled{%
  \def\@fs@cfont{\bfseries}\let\@fs@capt\floatc@ruled
  \def\@fs@pre{\vspace*{5pt}\hrule height.8pt depth0pt \kern2pt}%
  \def\@fs@post{\kern2pt\hrule\relax}%
  \def\@fs@mid{\kern2pt\hrule\kern2pt}%
  \let\@fs@iftopcapt\iftrue}
\floatstyle{betterruled}
\restylefloat{algorithm}
\makeatother
\usepackage[letterpaper, left=0.75in, right=0.75in, bottom=0.75in, top=0.75in]{geometry}

\title{\LARGE \bf
Inner and Outer Approximations of Star-Convex Semialgebraic Sets
}
\author{James Guthrie\thanks{J. Guthrie is with the Department of Electrical and Computer Engineering, 3400 N. Charles Street, Johns Hopkins University, Baltimore, MD 21218, USA. Email: {\tt\small jguthri6@jhu.edu} }}
\begin{document}

\maketitle
\thispagestyle{empty}
\pagestyle{empty}

%%%%%%%%%%%%%%%%%%%%%%%%%%%%%%%%%%%%%%%%%%%%%%%%%%%%%%%%%%%%%%%%%%%%%%%%%%%%%%%%
\begin{abstract}
We consider the problem of approximating a semialgebraic set with a sublevel-set of a polynomial function. In this setting, it is standard to seek a minimum volume outer approximation and/or maximum volume inner approximation. As there is no known relationship between the coefficients of an arbitrary polynomial and the volume of its sublevel sets, previous works have proposed heuristics based on the determinant and trace objectives commonly used in ellipsoidal fitting. For the case of star-convex semialgebraic sets, we propose a novel objective which yields both an outer and an inner approximation while minimizing the ratio of their respective volumes. This objective is scale-invariant and easily interpreted. Numerical examples demonstrate that the approximations obtained are often tighter than those returned by existing heuristics. We also provide methods for establishing the star-convexity of a semialgebraic set by finding inner and outer approximations of its kernel.
\end{abstract}

%%%%%%%%%%%%%%%%%%%%%%%%%%%%%%%%%%%%%%%%%%%%%%%%%%%%%%%%%%%%%%%%%%%%%%%%%%%%%%%%

\section{INTRODUCTION}
Consider a compact, semialgebraic set $\mathcal{X} \subset \mathbb{R}^n$ given by the intersection of the 1-sublevel sets of $m$ polynomial functions $g_i(x) \in \mathbb{R}[x]$:
\begin{equation} \label{eqn:setX}
    \mathcal{X} = \{x \,|\, g_i(x) \leq 1,\, i \in [m]\}.
\end{equation}
Semialgebraic sets arise naturally in many control applications. The set of coefficients for which a polynomial is Schur or Hurwitz stable is given by a semialgebraic set. \textcolor{\editcolor}{For Hurwitz stability, the polynomial inequalities can be derived from the Routh array.}
% In motion planning, a vehicle's occupied space $\mathcal{V}(q) \subset \mathbb{R}^n$ is a function of its configuration $q \in \mathbb{R}^x$. Collision avoidance with an obstacle $\mathcal{O} \subset \mathbb{R}^n$ can be posed as ensuring $0 \not\in \mathcal{M}$ where $\mathcal{M} = \{v - o \, | \, v \in \mathcal{V}(q), o \in \mathcal{O} \}$. The latter is a semialgebraic set arising from the Minkowski difference.
These sets are often complicated and cumbersome to analyze. As such, it is common to seek simpler representations which closely approximate the set but are more amenable to further analysis \cite{Dabbene2017}. Examples of ``simple'' representations include hyperrectangles and ellipsoids. 

A number of publications have explored the use of sum-of-squares (SOS) optimization for approximating a semialgebraic set with a simpler representation \cite{Dabbene2017, Magnani2005, Henrion2012, Ahmadi2017, Cerone2012, Guthrie2021, Guthrie2022, Jones2019}. The most common parameterization is to seek a SOS polynomial whose 1-sublevel set $\mathcal{F} = \{ x \, | \, f(x) \leq 1 \}$ provides either an inner ($\mathcal{F} \subseteq \mathcal{X}$) or outer ($\mathcal{F} \supseteq \mathcal{X}$) approximation of the set $\mathcal{X}$. 
\begin{figure}
    \centering
    \includegraphics[width=0.48\textwidth]{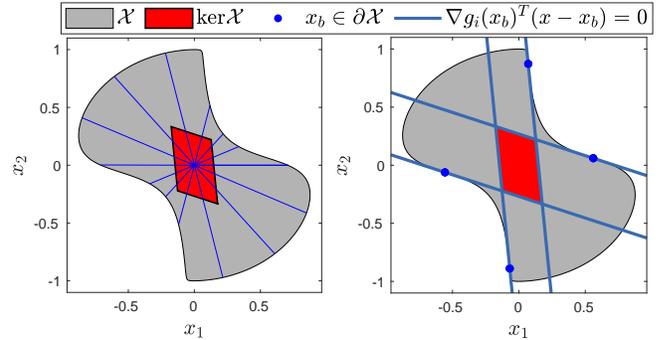}
    \caption{\textcolor{\editcolor}{The kernel is the convex set of points $p \in \mathcal{X}$ such that the line segment $\overline{pq} \subseteq \mathcal{X}$ for any $q \in \mathcal{X}$ (left)}. It is given by the intersection of all linearized active constraints $g_i(x_b) = 1$ defining $\partial\mathcal{X}$ (right).}
    \label{fig:star_convex}
    \vspace{-1mm}
\end{figure}
In this formulation, an open question is the choice of the objective function. For outer (resp. inner) approximations, a natural objective is to minimize (resp. maximize) the volume of the 1-sublevel set. For an ellipsoid $\mathcal{E} = \{x \, | \, x^TAx + b^Tx + c \leq 1 \}$ where $A \succeq 0$, the volume is proportional to $\textup{det} A^{-1}$. Using the logarithmic transform, ellipsoidal volume minimization can be posed as the convex objective $-\textup{logdet}A$ \cite{Magnani2005}. More generally, in the case of homogeneous polynomials it is possible to find the minimum volume outer approximation by solving a hierarchy of semidefinite programs \cite{Lasserre2014}. 
 
% While ellipsoids and homogeneous polynomials offer established techniques for approximating a set, they have inherent symmetry. Thus they are not ideal candidates for approximating general, asymmetric shapes.
Ellipsoids and homogeneous polynomials are not ideal candidates for approximating asymmetric shapes due to their inherent symmetry.
General polynomials offer a more flexible basis for approximating sets. The caveat is that we lack expressions for computing the volume of the 1-sublevel set as a function of the polynomial coefficients. The most common approach is to mimic the determinant (\cite{Magnani2005, Ahmadi2017}) or trace \cite{Dabbene2017} objectives used in ellipsoidal fitting. These objectives often yield qualitatively good approximations. 
However, they have no explicit relationship to the volume beyond upper bounding it in some cases \cite{Dabbene2017}. Thus it is difficult to infer the quality of an approximation from the objective value attained. 
% Instead, assessing the quality of the solution requires post-processing by either 1) numerically computing the resulting volume or 2) plotting the resulting set.

% They suffer from a lack of interpretability as the objective generally has no explicit relationship to the volume of the 1-sublevel set beyond upper bounding it in some cases \cite{Dabbene2017}. Assessing the quality of the solution requires post-processing by either 1) numerically computing the resulting volume or 2) plotting the resulting set. % for qualitative assessment.

\subsection{Contributions}
% We propose a new approach for finding inner and outer approximations of semialgebraic sets using SOS optimization. Our method is tailored to cases in which the set is star-convex. 
% To our knowledge, star-convexity has not been explored in the SOS literature with the exception of \cite{Wang2005}. 
% Our contributions are as follows:
\textcolor{\editcolor}{This paper makes the following contributions:}
\begin{itemize}
    \item We propose and justify an algorithm based on SOS optimization for jointly finding an inner and outer approximation of a semialgebraic set. The algorithm minimizes the volume of the outer approximation relative to the volume of the inner approximation. 
    \textcolor{\editcolor}{This objective is easily interpreted and scale-invariant.}
    \item
    We provide numerical examples showing that our algorithm tends to yield better approximations than existing methods when applied to star-convex sets.
    \item We provide algorithms for finding inner and outer approximations of the kernel of a star-convex set as shown in Figure \ref{fig:star_convex}.
\end{itemize}

The paper is organized as follows. Section II defines the problem we address and reviews the notion of star-convexity. Section III surveys existing volume heuristics for SOS-based set approximation. Section IV proposes a new volume heuristic for finding outer and inner approximations. Section V provides methods for approximating the kernel of a star-convex set. Section VI provides numerical examples. Section VII concludes the paper.

\newcommand{\Point}[1]{r_{#1}e^{i\theta_{#1}}}
\newcommand{\CPoint}[1]{\mathbf{r}_{#1}}
\newcommand{\Zero}{{0}}
\newcommand{\hull}{\textup{conv}}
\newcommand{\sos}[1]{\Sigma[#1] }
\newcommand{\obs}{\mathcal{O}}

\newcommand{\TrajOpt}{(14)}
\newcommand{\TrajOptDyn}{(14c)}

\newcommand{\Set}{\mathcal{S}}
\newcommand{\Ker}{\textup{ker}}
\newcommand{\A}{\mathcal{A}}
\newcommand{\B}{\mathcal{B}}

\newcommand{\F}{\mathcal{F}_1}
\newcommand{\Fk}{\mathcal{F}_{1,\kappa}}
\newcommand{\dX}{\partial\mathcal{X}}

\newcommand{\KerOuter}{\mathcal{K}_o}
\newcommand{\KerInner}{\mathcal{K}_i}
\newcommand{\SetInner}{\underaccent{\bar}{\mathcal{F}}_1}
\newcommand{\SetOuter}{\Tilde{\mathcal{F}}_1}
\newcommand{\scale}{s}
\subsection{Notation}
Let $i \in [k]:= \{1,\hdots,k\}$. Let $\mathbb{Z}^+$ denote the set of positive integers. Let $S^{n-1} := \{ x \in \mathbb{R}^n \, | \, \|x\| = 1 \}$. The notation $P \succeq 0$ indicates that the symmetric matrix $P$ is positive semidefinite (PSD). Given a compact set $\mathcal{X} \subset \mathbb{R}^n$, its volume (formally, Lebesgue measure) is denoted $\textup{vol } \mathcal{X} := \int_\mathcal{X} dx$. Let $\sigma_{\mathcal{X}}(c) := \underset{x \in \mathcal{X}}{\textup{max}} \, c^Tx$ denote the support function of $\mathcal{X}$ where $c \in S^{n-1}$. Given sets $\A,\B \subseteq \mathbb{R}^n$ the (bi-directional) Hausdorff distance is $d_H(\A,\B) := \textup{max}(h(\A,\B), h(\B,\A))$ where $h(\A,\B) := \underset{a \in \A}{\textup {max}} \, \underset{b \in \B}{\textup {min}} \, \|a-b\|_2 $. 

The $\alpha$-sublevel set of a function $f(x): \mathbb{R}^n \rightarrow \mathbb{R}$ is $\{ x \in \mathbb{R}^n \, | \, f(x) \leq \alpha \}$. For $x \in \mathbb{R}^n$, let $\mathbb{R}[x]$ denote the set of polynomials in $x$ with real coefficients.  Let $\mathbb{R}_d[x]$ denote the set of all polynomials in $\mathbb{R}[x]$ of degree less than or equal to $d$. A polynomial $p(x) \in \mathbb{R}[x]$ is a SOS polynomial if there exists polynomials $q_i(x) \in \mathbb{R}[x], i \in [j]$ such that $p(x) = q_1^2(x) + \hdots + q_j^2(x)$. We use $\sos{x}$ to denote the set of SOS polynomials in $x$.  A polynomial of degree $2d$ is a SOS polynomial if and only if there exists $P \succeq 0$ (the Gram matrix) such that $p(x) = z(x)^TPz(x)$ where $z(x)$ is the vector of all monomials of $x$ up to degree $d$ \cite{Parillo2000}. Letting $m := \binom{n+d}{d}$ denote the length of $z(x)$, we have that $P \in \mathbb{R}^{m \times m}$. To minimize notational clutter, we will sometimes list a polynomial $f(x)$ as a decision variable. It is implied that a degree is specified and matrix $P$ is introduced as a decision variable such that $f(x) = z(x)^TPz(x)$. 
\section{Problem Statement}
\begin{defn}[Star-Convex Set \textcolor{\editcolor}{\cite{Brunn1913}}]
A set $\mathcal{S} \subseteq \mathbb{R}^n$ is star-convex if it has a non-empty kernel. The kernel is
\begin{equation}
    \Ker \; \Set := \{ x \,|\, tx + (1-t)y \in \Set \,\forall\, t \in [0,1], y \in \Set \}.
\end{equation}
\end{defn}
The kernel is the set of points in $\mathcal{S}$ from which one can ``see'' all of $\mathcal{S}$ as shown in Figure \ref{fig:star_convex}. It is easily shown that the kernel is convex. If $\mathcal{S}$ is convex then $\Ker\, \mathcal{S} = \mathcal{S}$. 

We will be interested in approximating the set \eqref{eqn:setX} for the case in which it is star-convex with respect to the origin. 
\begin{problem}[Star-Convex Set Approximation]\label{prob:approx_star_cvx_outer_inner}
Given a compact, semialgebraic set $\mathcal{X}$ with $0 \in \textup{int}\mathcal{X} \cap \Ker\mathcal{X}$ and $d \in \mathbb{Z}^+$ find a polynomial $f_o(x) \in \mathbb{R}_{2d}[x]$ ($f_i(x) \in \mathbb{R}_{2d}[x]$) whose 1-sublevel set $\mathcal{F}_o$ ($\mathcal{F}_i$) is of minimum (maximum) volume and is an outer (inner) approximation of $\mathcal{X}$:
\begin{equation*}
% \underset{f_o(x) \in \mathbb{R}_{2d}[x]}{\textup{min}}  \textup{ vol } \mathcal{F}_o
% \textup{ s.t. } \mathcal{X} \subseteq \mathcal{F}_o 
\underset{f_o(x) \in \mathbb{R}_{2d}[x]}{\textup{min}}  \textup{ vol } \mathcal{F}_o
\textup{ s.t. } \mathcal{X} \subseteq \mathcal{F}_o 
\end{equation*}
\begin{equation*}
\left(
\underset{f_i(x) \in \mathbb{R}_{2d}[x]}{\textup{max}} \textup{ vol } \mathcal{F}_i
\textup{ s.t. } \mathcal{F}_i \subseteq \mathcal{X}
\right).
\end{equation*}

\end{problem}

% \begin{problem}[Minimum Volume Outer Approximation]\label{prob:approx_star_cvx_outer}
% Given a compact, semialgebraic set $\mathcal{X}$ with $0 \in \textup{int}\mathcal{X} \cap \Ker\mathcal{X}$ and $d \in \mathbb{Z}^+$ find a polynomial $f_o(x) \in \mathbb{R}_d[x]$ whose 1-sublevel set $\mathcal{F}_o$ is of minimum volume and contains $\mathcal{X}$.
% \begin{equation}
% \underset{f_o(x) \in \mathbb{R}_d[x]}{\textup{min}} \textup{ vol } \mathcal{F}_o
% \textup{ s.t. } \mathcal{X} \subseteq \mathcal{F}_o
% \end{equation}
% \end{problem}

% \begin{problem}[Maximum Volume Inner Approximation]\label{prob:approx_star_cvx_inner}
% Given a compact, semialgebraic set $\mathcal{X}$ with $0 \in \textup{int}\mathcal{X} \cap \Ker\mathcal{X}$ and $d \in \mathbb{Z}^+$ find a polynomial $f_i(x) \in \mathbb{R}_d[x]$ whose 1-sublevel set $\mathcal{F}_i$ is of maximum volume and is contained in $\mathcal{X}$.
% \begin{equation}
% \underset{f_i(x) \in \mathbb{R}_d[x]}{\textup{max}} \textup{ vol } \mathcal{F}_i
% \textup{ s.t. } \mathcal{F}_i \subseteq \mathcal{X}
% \end{equation}
% \end{problem}

To establish star-convexity of $\mathcal{X}$, we seek polytopic approximations of its kernel. \begin{problem}[Kernel Approximation]\label{prob:approx_kernel}
Given a semialgebraic set $\mathcal{X} \subset \mathbb{R}^n$ find a polytope $\KerOuter$ ($\KerInner$) of minimum (maximum) volume that is an outer (inner) approximation of $\Ker\mathcal{X}$:
\begin{equation*}
\textup{min vol }\KerOuter
\textup{ s.t. } \Ker\mathcal{X} \subseteq \KerOuter
\end{equation*}
\begin{equation*}
\left( \textup{max vol }\KerInner
\textup{ s.t. } \KerInner \subseteq \Ker\mathcal{X} \right).
\end{equation*}
\end{problem}

\section{Existing Volume Heuristics for Set Approximation}
We review existing heuristics for approximating semialgebraic set $\mathcal{X}$ using SOS optimization. Each of these methods finds an even-degree polynomial $f(x) = z(x)^TPz(x)$. The variations between the methods largely relate to the objective applied to Gram matrix $P$. For general polynomials, there is no known relationship between $P$ and the volume of the sublevel sets. Thus the following objectives are all heuristics in some sense.
% Each of these methods finds an even-degree polynomial $f(x) \in \mathbb{R}_{2d}[x]$ where $d \in \mathbb{Z}^+$ is specified by the user. The polynomial is parameterized as $f(x) = z(x)^TPz(x)$, where $z(x)$ is the vector of all monomials of $x$ up to degree $d$ and $P$ is the Gram matrix.  
% Here $z(x)$ is the vector of all monomials of $x$ up to degree
% $d$. is a monomial basis with $m$ terms and $P \in \mathbb{R}^{m \times m}$ is a symmetric matrix decision variable.   

\subsection{Determinant Maximization $(-\textup{det} P)$}
% For ellipsoids $\mathcal{E} = \{ x \, | \, x^TPx \leq 1 \}$ where $P \succ 0$, the volume of a sublevel set is proportional to $\textup{det} P^{-1}$. Thus given a set $\mathcal{X}$, one can find the minimum volume outer ellipsoid by maximizing the determinant of its Gram matrix $P$ subject to a set containment constraint $\mathcal{X} \subseteq \mathcal{E}$. Determinant maximization can be transformed into a convex objective through the logarithmic transform \cite{Boyd2004}.  

In \cite{Magnani2005}, the authors propose maximizing the determinant of the Hessian $\nabla^2f(x)$ of SOS polynomials. If $f$ is a polynomial of degree 2, this reduces to the ellipsoidal objective $-\textup{det} A$ for $\mathcal{E} = \{x \, | \, x^TAx + b^Tx + c \leq 1 \}, A \succeq 0$. As the Hessian must be PSD, the outer approximation is convex. This makes it ill-suited to approximating non-convex shapes. 

In \cite{Ahmadi2017}, the authors propose performing determinant maximization directly on the Gram matrix $P$. The Hessian is no longer required to be PSD. This allows non-convex outer approximations to be found. 

% \begin{mini}|s|
% {P}{-\textup{log det} P} {}{}
% \addConstraint{f(x) = z(x)^TPz(x), \;}{P \succeq 0,}{}
% \addConstraint{f(x) \leq}{1 \, \forall \, x \in \mathcal{X}} {}
% \label{} 
% \end{mini}

\subsection{Inverse Trace Minimization $(\textup{tr}P^{-1})$}
The determinant maximization objective minimizes the product of the eigenvalues of $P^{-1}$. In \cite{Ahmadi2017}, the authors propose an alternative heuristic of minimizing the sum of the eigenvalues of $P^{-1}$. This requires an additional matrix variable $V$ and constraint $V \succeq P^{-1}$. Using the Schur complement this can be written as a block matrix constraint involving $V$ and $P$ (vice $P^{-1}$). The objective $\textup{min tr}V$ then indirectly minimizes the sum of the eigenvalues of $P^{-1}$.

% \begin{mini}|s|
% {P,V}{-\textup{log det} P} {}{}
% \addConstraint{f(x) = z(x)^TPz(x), \;}{\begin{bmatrix} V & I \\ I & P \end{bmatrix} \succ 0,}{}
% \addConstraint{f(x) \leq}{1 \, \forall \, x \in \mathcal{X}} {}
% \label{} 
% \end{mini}

\subsection{$l_1$ Minimization}
In \cite{Dabbene2017} the authors propose minimizing the $l_1$ norm of a polynomial evaluated over a bounding box $\mathcal{B} \supseteq \mathcal{X}$. This approach was first introduced in~\cite{Henrion2009} for approximating the volume of semialgebraic sets.
Using hyperrectangles as bounding boxes, one can integrate the polynomial over $\mathcal{B}$. The resulting objective $l_1(f(x)) := \int_{\mathcal{B}} f(x) \, dx$ is linear in terms of $P$. 
%\textcolor{\editcolor}{and provides an upper bound on the volume of the approximation}. 
The outer approximation consists of the intersection of the 1-superlevel set of $f(x)$ and $\mathcal{B}$:
\begin{equation}
    \mathcal{X} \subseteq ( \mathcal{B} \cap \{ x \,|\, f(x) \geq 1 \} ).
\end{equation}
This differs from other objectives which do not rely on bounding boxes as part of the set approximation.\footnote{One application of approximating semialgebraic sets is to yield a single sufficient condition for ensuring $x \not\in \mathcal{X}$, which can be incorporated into a nonlinear optimization problem (e.g. obstacle avoidance in motion planning \cite{Guthrie2022}). The presence of the bounding box in the resulting set description would require logical constraints to represent ($f(x) < 1 \lor x \not\in \mathcal{B}) \implies x \not\in \mathcal{X}$ which are generally unsupported in nonlinear optimization solvers.} In this setting, $f(x)$ is approximating the indicator function of $\mathcal{X}$ over a compact set $\mathcal{B}$. Convergence of $f(x)$ to the true indicator function in the limit (as degree $d \rightarrow \infty $) can be shown by leveraging the Stone-Weierstrass theorem. \textcolor{\editcolor}{The asymptotic rate of convergence is at least $O(1 / \textup{log log } d)$ \cite{Korda2018}}.
% However, the rate of convergence remains unknown and for a given degree $d$, other objectives may return a better outer approximation. 
Inner approximations can be found by outer approximating the complement of $\mathcal{X}$.
% \begin{mini}|s|
% {P}{\int_{\mathcal{B}} f(x) \, dx} {}{}
% \addConstraint{f(x) \geq}{0}{\forall x \in \mathcal{B},}
% \addConstraint{f(x) \geq}{1}{\forall x \in \mathcal{X}}
% \label{} 
% \end{mini}

\section{Inner and Outer Approximations of Star-Convex Sets}
\begin{figure}
    \centering
    \vspace{2mm}
    \includegraphics[width=0.48\textwidth, trim={0 0 0 0},clip]{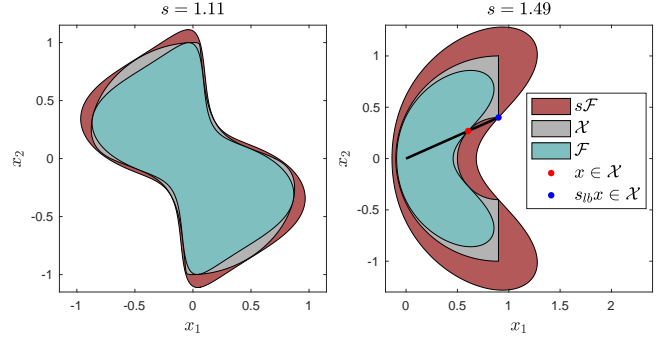}
    \vspace{-0.75cm}
    \caption{4th-order approximations of star-convex set (left) and non-star-convex set (right) found by minimizing scaling term $s$. The non-star-convex set has a lower bound $s_{lb} > 1$ on the achievable approximation scaling $s$.}
    \label{fig:scale_example}
    \vspace{-0.25cm}
\end{figure}

 We propose a new volume heuristic for solving Problem \ref{prob:approx_star_cvx_outer_inner}. Our heuristic is inspired by the following two lemmas.
 %The following lemma can be shown via a change of variables.
 \begin{lem}
Let $\mathcal{X}, \mathcal{F}$ be compact sets in $\mathbb{R}^n$ such that $\mathcal{F} \subseteq \mathcal{X}$. Let $0 \in \textup{int } \mathcal{F}$. Then there exists a scaling $\scale \geq 1$ such that $\mathcal{X} \subseteq \scale\mathcal{F}$.
\end{lem}
\begin{lem}\label{lem:volume_scale}
Let $\mathcal{X} \subset \mathbb{R}^n$. Let $\scale\mathcal{X} = \{ sx \,|\, x \in \mathcal{X} \}$ denote the scaled set where $\scale \geq 0$. Then $\textup{vol } \scale\mathcal{X} = s^n\cdot\textup{vol}\mathcal{X}$.
\end{lem}
% The following lemma provides conditions under which we can scale an inner approximation to become an outer approximation of a given compact set.
Thus given an inner approximation $\mathcal{F}$, we can obtain an outer approximation $\scale \mathcal{F}$ for some $s \geq 1$ with relation

%Lemma \ref{lem:volume_scale} gives the following relation:
%From Lemma \ref{lem:volume_scale}, we have the following relation:
\begin{equation}
    \frac{\textup{vol } \scale \mathcal{F}} {\textup{vol } \mathcal{F}} = \scale^n.
\end{equation}
By minimizing $\scale$ we minimize the ratio of the outer approximation volume to the inner approximation volume. Figure \ref{fig:scale_example} visualizes this intuitive heuristic for approximating a set.
% Taken together, these lemmas suggest an intuitive heuristic for jointly finding an inner and outer approximation of a compact set $\mathcal{X}$ by minimizing the scaling required to turn an inner approximation into an outer approximation. 

% INNER APPROXIMATION
We seek a polynomial $f : \mathbb{R}^n \rightarrow \mathbb{R}$ whose 1-sublevel set $\mathcal{F} = \{ x \, | \, f(x) \leq 1 \}$ is an inner approximation of $\mathcal{X}$. We turn this into a condition involving the complement of $\mathcal{X}$:
\begin{equation}
    \mathcal{F} \subseteq \mathcal{X} \iff f(x) > 1 \,\forall\, x \in \mathcal{X}^{\mathsf{c}}. \label{eqn:inner_comp_condition}
\end{equation}
Optimization methods require non-strict inequalities. We approximate the strict inequality by introducing a small constant $\epsilon > 0$ and working with the closure of the complement of $\mathcal{X}$. Define the following:
\begin{equation}
    \bar{\mathcal{X}} = \underset{i \in [m]}{\bigcup} \{ x \,|\, g_i(x) \geq 1 \}.
\end{equation}
We then use the following approximation of \eqref{eqn:inner_comp_condition}:
\begin{equation}
    \mathcal{F} \subset \textup{int}\mathcal{X} \Leftarrow f(x) \geq 1 + \epsilon \; \forall \, x \in \bar{\mathcal{X}}.
\end{equation}

Next, we scale the set $\mathcal{F}$ by a scaling variable $s > 1$ to obtain an outer approximation:
\begin{equation}
    s\mathcal{F} \supseteq \mathcal{X} \iff f(\tfrac{x}{s}) \leq 1 \,\forall\, x \in \mathcal{X}.
\end{equation}

Combining the above we arrive at the following:
\begin{mini}|s|
{f(x), \scale}{\scale} {}{}
\addConstraint{f(x)}{\geq 1+\epsilon}{\; \forall \, x \in \bar{\mathcal{X}},}
\addConstraint{f(\tfrac{x}{s})}       {\leq 1}         {\; \forall \, x \in \mathcal{X}.}
\label{opti:min_s} 
\end{mini}

% \begin{figure}
%     \centering
%     \includegraphics[width=0.5\textwidth, trim={0 0 0 0},clip]{scale_example.eps}
%     %\setlength{\belowcaptionskip}{-28pt} 
%     \caption{4th-Order Approximations of Star-Convex Set}
%     \label{fig:scale_example}
% \end{figure}

%For example, in Figure \ref{fig:scale_example} we have $\frac{2.0183}{1.6307} = 1.1125^2$
\begin{rem}
Our scaling heuristic is applicable to approximating any compact set containing the origin in its interior. However, it is best suited to approximating star-convex sets in which $0 \in \textup{int} \mathcal{X} \cap \Ker \, \mathcal{X} $ as visualized in Figure \ref{fig:scale_example}. \textcolor{\editcolor}{Otherwise there exists a lower bound $s_{lb}$ such that $1 < s_{lb} \leq s$ in \eqref{opti:min_s}.}
\begin{lem}\label{lem:s_lowerbound}
\textcolor{\editcolor}{Let $\mathcal{X}$ and $\mathcal{F}$ be compact sets in $\mathbb{R}^n$. Let $\mathcal{F} \subseteq \mathcal{X} \subseteq s^\star \mathcal{F}$ for some $s^\star > 1$. Let $0 \in \textup{int } \mathcal{F}$. 
Let $x, sx \in \mathcal{X}$ and $tx \not\in \mathcal{X} \, \forall \, t \in (1,s)$ for some $s > 1, x \neq 0$. Then $s^\star \geq s$. }
% \textcolor{\editcolor}{Let $\mathcal{X}, \mathcal{F}$ be compact sets in $\mathbb{R}^n$ such that $\mathcal{F} \subseteq \mathcal{X}$. Let $0 \in \textup{int } \mathcal{F}$. Let $\mathcal{X} \subseteq s\mathcal{F}$ for some $s > 1$. Then it must hold that $s \geq s_{lb}$ where $s_{lb} = \textup{sup}\{\frac{1}{a} \, | \, 0 < a < 1, x \neq 0, ax, x \in \mathcal{X}, tx \not\in \mathcal{X} \, \forall \, t \in (a, 1) \}$.} 
\begin{proof}
% See \cite{Guthrie2022_star_arxiv}.\phantom\qedhere
See appendix.\phantom\qedhere
\end{proof}
\end{lem}
\textcolor{\editcolor}{
We let $s_{lb}$ denote the greatest lower bound given by Lemma \ref{lem:s_lowerbound}. This imposes a minimum volume ratio between the inner and outer approximation. Figure \ref{fig:scale_example} (right) visualizes this result. The set is not star-convex and therefore $0 \not \in \Ker \mathcal{X}$. The black line segment connecting the origin to point $s_{lb}x$ is not contained in $\mathcal{X}$. This point imposes a lower bound on $s$, preventing the inner and outer approximations from coming closer together.
}
\end{rem}
We introduce SOS polynomials $\lambda_i(x), \mu_i(x), i \in [m]$ and replace the set-containment conditions in \eqref{opti:min_s} with SOS conditions.\footnote{\textcolor{\editcolor}{For the outer approximation of the compact set $\mathcal{X}$, the SOS conditions are necessary and sufficient by Putinar's Positivstellensatz when $\mu_i(x)$ is of high-enough degree and the defining polynomials $g_i$ satisfy the Archimedean assumption \cite{Putinar1993}. The inner approximation constraint involves an unbounded set. The associated SOS reformulation utilizes the generalized $\mathcal{S}$-procedure which is only sufficient \cite{Parillo2000}. 
% Thus for a fixed $s$, \eqref{opti:min_s} may be feasible but $\textup{FindApprox}(s,g_i)$ infeasible, making our approach suboptimal.
}}
If $\scale$ is left as a decision variable, we would have bilinear terms involving the coefficients of $f(x)$ and $\scale$. Instead we perform a bisection over $\scale$, solving a feasibility problem at each iteration as given by \eqref{opti:FindApprox}. 
Algorithm \ref{alg:inner_outer_X} details the bisection method.

\noindent\textbf{Optimization Problem: FindApprox($\scale, g_i$)}
\begin{mini}|s|
{f(x), \lambda_i(x), \mu_i(x)}{0} {}{}
\addConstraint{f(x) - (1+\epsilon) - \lambda_i(x)(g_i(x)-1)}{\in \sos{x},}{\,i \in [m],}
\addConstraint{1-f(\tfrac{x}{s}) - \sum_{i=1}^m \mu_i(x)(1-g_i(x))}{\in \sos{x},}{}
\addConstraint{\lambda_i(x), \mu_i(x)}{\in \sos{x}, \; \; \; \; \; i \in [m].}
\label{opti:FindApprox} 
\end{mini}

 \begin{algorithm} 
 \caption{Inner and Outer Approximation of $\mathcal{X}$}
 \begin{algorithmic}[H]\label{alg:inner_outer_X}
 \renewcommand{\algorithmicrequire}{\textbf{Input:}}
 \renewcommand{\algorithmicensure}{\textbf{Output:}}
 %\REQUIRE $\mathcal{X} \subset \mathbb{R}^n$, $z(x) \in \mathbb{R}[x]$, $s_{tol} > 0$
 \REQUIRE $\mathcal{X} = \{x \in \mathbb{R}^n \,|\, g_i(x) \leq 1, i \in [m]\}$, $s_{tol} > 0$
 %\REQUIRE $g_i(x), i \in [m]$, $s_{tol} > 0$
 \ENSURE  $\mathcal{F}, s\mathcal{F} \textup{ s.t. } \mathcal{F} \subseteq \mathcal{X} \subseteq s\mathcal{F}$
  \STATE $s_{ub} \gets 1+s_{tol}, s_{lb} \gets 1$
  \WHILE {FindApprox($s_{ub}, g_i) = \textup{Infeasible}$}
  \STATE $s_{lb} \gets s_{ub}$
  \STATE $s_{ub} \gets 2 s_{ub}$
  \ENDWHILE
  \WHILE{$s_{ub} - s_{lb} > s_{tol}$}
  \STATE $s_{try} \gets 0.5(s_{ub} + s_{lb})$
  \IF{FindApprox($s_{try}, g_i) = \textup{Infeasible}$}
  \STATE $s_{lb} \gets s_{try}$
  \ELSE
  \STATE $s_{ub} \gets s_{try}$
  \ENDIF
  \ENDWHILE
 \RETURN FindApprox($s_{ub}, g_i)$
 \end{algorithmic} 
 \end{algorithm}
\vspace{-0.25cm}

\begin{rem}
The objective is scale-invariant. Let solution $(f^*(x), \scale^*)$ define an outer and inner approximation of $\mathcal{X}$. Scale $\mathcal{X}$ by $\alpha > 0$, 
\textcolor{\editcolor}{replacing constraints $g_i(x)$ with $g_i(\frac{x}{\alpha})$.} 
Then the solution pair $(f^*(\tfrac{x}{\alpha}), \scale^*)$ defines the new approximation, where the objective value remains unchanged. \textcolor{\editcolor}{The objective is not translation-invariant however. For example, assume we approximate a star-convex set exactly with $(f^\star(x), s^\star = 1)$. Translate $\mathcal{X}$ by $t \in \mathcal{X} \setminus \Ker \mathcal{X}$, replacing $g_i(x)$ with $g_i(x-t)$. Then $0 \not \in \Ker \mathcal{X}$ and $s^\star > 1$ for any approximation by Lemma \ref{lem:s_lowerbound}.}
\end{rem}
% This is shown by the following relation which is trivial to verify.
% \begin{equation}
% \mathcal{F} \subseteq \mathcal{X} \subseteq \scale\mathcal{F} \Leftrightarrow
% \alpha \mathcal{F} \subseteq \alpha \mathcal{X} \subseteq \alpha \scale \mathcal{F}
% \end{equation}

\begin{rem}
If $\mathcal{F}$ is convex we can relate the scaling $s$ to the Hausdorff distance between the approximations.
\begin{lem}\label{lem:hausdorff}
Let $\mathcal{F} \subset \mathbb{R}^n$ be a convex, compact set and $s \geq 1$. Then the following holds: 
\begin{equation}
    %d_H(s\mathcal{F},\mathcal{F}) = 
    % \underset{x \in s\mathcal{F}}{\textup {max}} \, \underset{y \in \mathcal{F}}{\textup {min}} \, \|x-y\|_2 = (s-1)\cdot \underset{x \in \mathcal{F}}{\textup {max}} \|x\|_2
        d_H(s\mathcal{F},\mathcal{F}) = (s-1)\cdot \underset{x \in \mathcal{F}}{\textup {max}} \|x\|_2.
\end{equation}
\end{lem}
\begin{proof}
% See \cite{Guthrie2022_star_arxiv}.\phantom\qedhere
See appendix.\phantom\qedhere
\end{proof}
\end{rem}

% \begin{rem}
% Existing approaches for finding inner or outer approximations with SOS optimization generally require $f(x)$ be SOS. Our formulation is slightly less restrictive as $f(x)$ can be negative for $x \in \mathcal{X}$. 
% %Therefore we search over a larger candidate set of polynomials. 
% We have found adding the constraint $f(x) \in \sos{x}$ does not noticeably degrade the quality of our approximations. \textcolor{\editcolor}{Without this constraint we can replace the Gram matrix representation $f(x) = z(x)^TPz(x)$ with a coefficient vector representation $f(x) = c^Tw(x)$. Here $w(x)$ is the vector of all monomials of $x$ up to degree $2d$. The coefficient vector $c$ contains $\binom{n+2d}{2d}$ terms.  The Gram matrix consists of $\frac{m(m+1)}{2}$ terms where $\binom{n+d}{d}$. The vector representation is more compact for $d > 1$. Further, while the Gram matrix is not unique \cite{Topcu2010}, the vector $c$ is unique.}
% \end{rem}

\section{Sampling-Based Approximations of the Kernel}
\textcolor{\editcolor}{
Algorithm \ref{alg:inner_outer_X} assumed the set $\mathcal{X}$ contained the origin in its kernel. If this does not hold, but there exists a point $x^\star \in \Ker \mathcal{X} \cap \textup{int}\mathcal{X}$ we can apply Algorithm \ref{alg:inner_outer_X} to the translated set $\{ x - x^* \, | \, x \in \mathcal{X} \}$. As our objective is not invariant with respect to translation, it is useful to approximate the kernel to establish possible choices for $x^\star$.\footnote{\textcolor{\editcolor}{A practical heuristic is to let $x^\star$ be the Chebyshev center of $\Ker \mathcal{X}$.}}
In this section we provide algorithms for finding polytopic approximations of $\Ker \, \mathcal{X}$.
}

It will be convenient to represent the boundary of $\mathcal{X}$ in terms of the inequality that is active. Define the following:
\begin{equation}
    \partial\mathcal{X}_i = \{x \,|\, g_i(x) = 1, g_j(x) \leq 1, j \in [m]\setminus i\}.
\end{equation}
The boundary of $\mathcal{X}$ is given by the union
\begin{equation}    
    \partial\mathcal{X} = \bigcup_{i \in [m]} \partial\mathcal{X}_i.
\end{equation}

\begin{lem}\label{lem:kernel_gradient_condition}
Let $\mathcal{X}$ be a semialgebraic set as defined in \eqref{eqn:setX}. \textcolor{\editcolor}{Let $\nabla g_i(x_b) \neq 0 \, \forall \, x_b \in \partial \mathcal{X}_i, i \in [m]$}. The kernel of $\mathcal{X}$ is given by the following semialgebraic set:
\begin{equation*}\label{eq:kernel}
\Ker\mathcal{X} = \{x_k \,|\, \nabla g_i(x_b)^T(x_k - x_b) \leq 0 \,\forall \, x_b \in \partial \mathcal{X}_i, i \in [m] \}.
\end{equation*}
\begin{proof}
%See \cite{Guthrie2022_star_arxiv}.\phantom\qedhere
See appendix.\phantom\qedhere
\end{proof}
\end{lem}
% $\Rightarrow:$
% Assume $x_k \in \Ker\mathcal{X}$ but there exists a point $x_b \in \partial\mathcal{X}_i$ for some $i \in [m]$ such that $\nabla g_i(x_b)^T(x_k - x_b) > 0$.
% Recall the definition of the directional derivative:
% \begin{equation*}
% \lim_{t\rightarrow0} \frac{g_i(tx_k + (1-t)x_b) - g_i(x_b)}{t} =  \nabla g_i(x_b)^T(x_k - x_b)    
% \end{equation*}
% Noting that $g_i(x_b) = 1$ and $\nabla g_i(x_b)^T(x_k - x_b) > 0$ yields
% \begin{equation*}
% \lim_{t\rightarrow0} \frac{g_i(tx_k + (1-t)x_b) - 1}{t} > 0    
% \end{equation*}
% This implies there exists an open interval $t \in (0,\alpha), \alpha > 0$ in which $g_i(tx_k + (1-t)x_b) > 1$. The associated line segment does not belong to $\mathcal{X}$, i.e. $\{tx_k + (1-t)x_b \, | \, t \in (0,\alpha) \} \not\subseteq \mathcal{X}$ and therefore $x_k \not\in\Ker\mathcal{X}$, contradicting our assumption.\newline
% $\Leftarrow$: The proof of the reverse direction is nearly identical.

\begin{remark}
From Lemma \ref{lem:kernel_gradient_condition} we see that the kernel of $\mathcal{X}$ is defined by cutting-planes tangent to the active constraint $g_i(x_b) = 1, x_b \in \partial\mathcal{X}$ as shown in Figure \ref{fig:star_convex}. 
\end{remark}
\begin{remark}
\textcolor{\editcolor}{Lemma \ref{lem:kernel_gradient_condition} assumes the gradient of an active constraint is non-zero. While restrictive, we note that this assumption is typically satisfied in sets of practical interest.}
\end{remark}
% In Section XXX we assume knowledge of a point $x^* \textbf{int}\mathcal{X}, x^* \in \Ker\mathcal{X}$ such that we can define a translated set $\mathcal{Z}$ containing the origin in its interior and kernel. The question arises as to 1) how to find the kernel of $\mathcal{X}$ and 2) how to pick a good point within the kernel to translate to the origin. We note that this also addresses the simpler question of testing whether $\mathcal{X}$ is star-convex as this corresponds to showing that $\Ker\mathcal{X} \neq \emptyset$.
%  Additionally, we provide a method for approximating the maximium inscribed ball. Empirically we have found that the center of this ball is a good point for defining the translated set $\mathcal{Z}$.

% Although $\Ker\mathcal{X}$ is a convex, semialgebraic set it is not straightforward to represent it within a semidefinite program. Determining if a convex, semialgebraic set is semidefinite-representable is an area of active research and there is not a systematic procedure for constructing the representation if one exists \cite{Helton2010}. Instead, 
We provide sampling-based algorithms for finding outer and inner approximations of this set. If the outer approximation is empty, this is sufficient to conclude that the set $\mathcal{X}$ is not star-convex. Conversely, if the inner approximation is not empty this is sufficient to establish that $\mathcal{X}$ is star-convex. In the case that the outer approximation is not empty and the inner approximation is empty we cannot conclude anything about the star-convexity of the set.

\subsection{Outer Approximation}
We assume the existence of an oracle $\text{Sample}(\partial\mathcal{X})$ which allows us to randomly sample points $x_b \in \partial\mathcal{X}$ and identify the set of active constraints $\mathcal{I} = \{i \, | \, i \subseteq [m], g_i(x_b) = 1 \}$.\footnote{\textcolor{\editcolor}{Starting from a point in the interior of $\mathcal{X}$, one can choose a direction and find a boundary point via bisection. Alternatively, nonlinear optimization methods may be leveraged to find boundary points.}}
From Lemma \ref{lem:kernel_gradient_condition}, each sample defines a cutting plane satisfied by $\Ker\mathcal{X}$. We collect these constraints to form an outer approximation $\KerOuter \supseteq \Ker\mathcal{X}$. If at any point, $\KerOuter = \emptyset$ (which can be determined using Farkas' Lemma) we terminate as this implies $\Ker\mathcal{X} = \emptyset$. Algorithm \ref{alg:outer_Ker} summarizes the method.

\vspace{-0.25cm}
 \begin{algorithm}
 \caption{Outer Approximation of $\Ker\mathcal{X}$}
 \begin{algorithmic}[H]\label{alg:outer_Ker}
 \renewcommand{\algorithmicrequire}{\textbf{Input:}}
 \renewcommand{\algorithmicensure}{\textbf{Output:}}
 \REQUIRE $\mathcal{X} = \{x \in \mathbb{R}^n \,|\, g_i(x) \leq 1, i \in [m]\}$, $n_s \geq 1$
 \ENSURE Outer Approximation $\KerOuter \supseteq \Ker\mathcal{X}$
  \STATE $\KerOuter \gets \mathbb{R}^n$
  \FOR {$j = 1$ to $n_s$}
  \STATE $x_b, \mathcal{I} \gets \text{Sample}(\partial\mathcal{X})$
  \STATE $\KerOuter \gets \KerOuter \bigcap \{ x \,|\, \nabla g_i^T(x_b)(x - x_b) \leq 0, i \in \mathcal{I} \}$
  \IF {($\KerOuter = \emptyset$)}
  \RETURN $\KerOuter$
  \ENDIF
  \ENDFOR
 \RETURN $\KerOuter$ 
 \end{algorithmic} 
 \end{algorithm}
\vspace{-0.25cm}

\subsection{Inner Approximation}
Consider finding a point $x_k \in \Ker \X$ that maximizes a linear cost $c^Tx_k$ where $c \in S^{n-1}$ (i.e. the support function of $\Ker\X$). From Lemma \ref{lem:kernel_gradient_condition}, the resulting convex optimization problem requires set containment constraints:
\begin{maxi}|s|
{x_k}{c^Tx_k} {}{}
\addConstraint{-\nabla g_i(x)^T(x_k - x) }{\geq 0 \, \forall \, x \in \partial\mathcal{X}_i, } {\, i \in [m]}.
\label{opti:kernel_inner} 
\end{maxi}

We replace the set containment conditions with SOS conditions \textcolor{\editcolor}{using Putinar's Positivstellensatz \cite{Putinar1993}.}
\newline
\noindent\textbf{Optimization Problem: FindSupport$(c, g_i)$}
\begin{maxi}|s|
{x_k, \lambda_j^{(i)}(x)}{c^Tx_k} {}{}
% \addConstraint{-\nabla g_i(x)^T(x_k - x) - \displaystyle\sum_{j=1}^m \lambda_j^{(i)}(x)(1-g_j(x))}{\in \sos{x},}{i \in [m]}
\addConstraint{-\!\nabla g_i(x)^T(x_k \! -   x) \! - \! \displaystyle\sum_{j=1}^m \lambda_j^{(i)}(x)(1 \!- \! g_j(x))}{\in \sos{x},}{\, i \in [m]}
\addConstraint{\lambda_j^{(i)}(x)}{\in \sos{x} , \; \; \, i \in [m], j \in [m] \setminus i}.
\label{opti:min_kappa} 
\end{maxi}

For a given direction $c \in S^{n-1}$ this program lower bounds the support function of $\Ker\mathcal{X}$. \textcolor{\editcolor}{The lower bound monotonically increases with $\textup{deg}(\lambda_j^{(i)})$.} If the problem is feasible, the maximizing argument $x_k$ belongs to $\Ker\mathcal{X}$ 
%(though it may lie in the interior) 
and therefore $\mathcal{X}$ is star-convex. If infeasible we cannot make any conclusions about the star-convexity of $\mathcal{X}$.
% \textcolor{\editcolor}{In this case we may increase the degree of $\lambda$ if we suspect the set is star convex.} 
By solving for random directions $c_i \in S^{n-1}, i \in [n_s]$ the convex hull of points $x_k$ provides an inner approximation of the kernel as given by Algorithm \ref{alg:inner_approx_ker}.%\footnote{\textcolor{red}{Under the assumption that $\partial \mathcal{X}$ is Archimedean, \eqref{opti:min_kappa} solves \eqref{opti:kernel_inner} for $\lambda_j$ of high enough degree.}}

\vspace{-0.25cm}
 \begin{algorithm}
 \caption{Inner Approximation of $\Ker\mathcal{X}$}\label{alg:inner_approx_ker}
 \begin{algorithmic}[H]\label{alg:inner_Ker}
 \renewcommand{\algorithmicrequire}{\textbf{Input:}}
 \renewcommand{\algorithmicensure}{\textbf{Output:}}
 \REQUIRE $\mathcal{X} \! = \! \{x \! \in \! \mathbb{R}^n | g_i(x) \! \leq \! 1, i \! \in \! [m]\}$,\! 
 $\{c_i\} \! \subset \! S^{n-1}, i \! \in \! [n_s]$
 \ENSURE Inner Approximation $\KerInner \subseteq \Ker\mathcal{X}$
  \STATE $\KerInner \gets \emptyset$
  \FOR {$j = 1$ to $n_s$}
  \STATE $x_k \gets \textup{FindSupport}(c_j, g_i)$
  \IF {$\textup{FindSupport}(c_j, g_i)$ = Infeasible}
  \RETURN $\KerInner = \emptyset$
  \ENDIF
  \STATE $\KerInner \gets \hull (\KerInner, x_k)$
  \ENDFOR
 \RETURN $\KerInner$ 
 \end{algorithmic} 
 \end{algorithm}
\vspace{-0.25cm}

% \begin{rem}
% If we first find an outer approximation $\KerOuter$ we can use $\Ker \mathcal{X} = \Ker \mathcal{X} \cap \KerOuter$
% we can use the following constraint:
% \begin{equation}
%     -\nabla g_i(x)^T(x_k - x) \geq 0 \forall x \in \partial\mathcal{G}_i^{\mathcal{X}} \cap \KerOuter, \, i \in [n]
% \end{equation}
% Adding these redundant constraints can often be advantageous when applying the $\mathcal{S}$-procedure.
% \end{rem}

\subsection{Kernel of Unions and Intersections}
Given sets $\mathcal{A}, \mathcal{B} \subseteq \mathbb{R}^n$ and their kernels, we can find inner approximations of the kernel of their intersection and union using the following lemma.
\begin{lem}\label{lem:kernel_relation}
Let $\mathcal{A}, \mathcal{B} \subseteq \mathbb{R}^n$. Then the following holds: 
% \footnote{Simple examples can be constructed to show that there is no relation between $\Ker (\mathcal{A} \cap \mathcal{B})$ and $\Ker ( \A \cup \B )$ in general.}
\begin{align}
\Ker ( \A \cap \B ) &\supseteq \Ker\A \cap \Ker\B \\
\Ker (\mathcal{A} \cup \mathcal{B}) &\supseteq \Ker\A \cap \Ker\B.
\end{align}
\end{lem}
\begin{proof}
% See \cite{Guthrie2022_star_arxiv}.\phantom\qedhere
See appendix.\phantom\qedhere
\end{proof}
Thus if $\mathcal{A}, \mathcal{B}$ are star-convex and have kernels that intersect, their union and intersection is also star-convex. This is useful for establishing star-convexity without resorting to numerical algorithms.

\section{Examples}
We evaluate Algorithm \ref{alg:inner_outer_X} on various examples and compare the results to the existing heuristics reviewed in Section III.\footnote{For the bounding box $\mathcal{B}$ required by the $l_1$ objective, we used the smallest hyperrectangle $\mathcal{B} \supseteq \mathcal{X}$ unless noted otherwise.} We focus our comparison on outer approximations as more heuristics apply to this case. We use percent error as our metric, calculated as $100 \times \frac{\textup{vol} \mathcal{F}_o - \textup{vol}\mathcal{X}} {\textup{vol}\mathcal{X}}$ where $\mathcal{F}_o$ is the outer approximation of $\mathcal{X}$. We first consider approximating two examples from the literature with polynomials of increasing degree. In all instances, our algorithm yielded the tightest outer approximation as shown in Figure \ref{fig:timing_tiled}.\footnote{We forego comparing 2nd-order polynomials as the determinant maximization objective exactly minimizes volume in this case.} Next we consider 100 randomly generated convex polytopes in $\mathbb{R}^2$. In the majority of cases, our heuristic yielded the tightest outer approximation as shown in Table \ref{tab:batch_cvx_2d}. \textcolor{\editcolor}{Lastly, we approximate a set that is not star-convex. Our heuristic degrades with increasing lower bound $s_{lb}$ as suggested by Lemma \ref{lem:s_lowerbound}.}
%Further, the $l_1$ objective did not yield improvements on the bounding box outer approximation $\mathcal{B}$ for this case.

% THIS IS THE TABLE I ORIGINALLY HAD
% \begin{table}
% \vspace{1.75mm}
% \caption{Percent Error of Outer Approximations of Examples A-C}
% \begin{center}
% \begin{tabular}{c c c c c c}  
% \hline
%   Example & Degree & $s$ & $-\textup{det} P$ & $\textup{tr} P^{-1}$ & $l_1$ \\ %[0.5ex] 
%   \hline
%   A & 4 & 11.9 & 35.1 & 40.0 & 18.3 \\
%   A & 6 & 1.4  & 8.3  & 10.0 & 12.8 \\ 
%   B & 4 & 17.7 & 31.1 & 35.0 & 37.3 \\
%   B & 6 & 4.9  & 9.7  & 14.0 & 17.7 \\ 
%   C & 4 & 2.6  & 20.1 & 21.2 & 15.3 \\
%   C & 6 & 0.6  & 7.2  &  7.4 & 11.0 \\
%  \hline
% \end{tabular}
% \label{tab:examples}
% \end{center}
% \vspace{-0.5cm}
% \end{table}

% \begin{table}
% \vspace{1.75mm}
% \caption{Percent Error of Outer Approximations of Examples A-C}
% \begin{center}
% \begin{tabular}{c c c c c}  
% \hline
%   Objective & Deg. & Ex. A & Ex. B & Ex. C \\ %[0.5ex] 
%   \hline
%  $s$ & 4 & 11.9 & 17.7  & 2.6 \\

%  $-\textup{det} P$  & 4 & 35.1 &  31.1 & 20.1 \\
 
%  $\textup{tr} P^{-1}$  & 4 & 40.0 &   35.0 &  21.2\\
 
%   $l_1$  & 4 & 18.3 &   37.3 &   15.3\\
%  \hline
%   $s$ & 6 & 1.4 &   4.9   &  0.6 \\

%  $-\textup{det} P$  & 6 & 8.3  &  9.7  &  7.2 \\

%  $\textup{tr} P^{-1}$  & 6 & 10.0  &  14.0 &   7.4\\

%   $l_1$  & 6 & 12.8 & 17.7  & 11.0\\
%  \hline
% \end{tabular}
% \label{tab:examples}
% \end{center}
% \vspace{-0.5cm}
% \end{table}

\subsection{Polynomial matrix inequality \textup{\cite{Henrion2012}}}
\begin{figure}
    \centering
    \vspace{0.2cm}
    \includegraphics[width=0.47\textwidth]{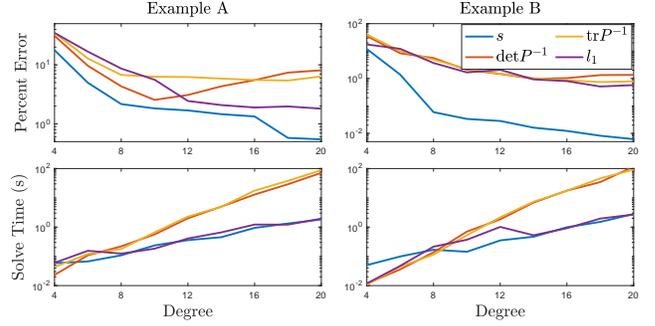}
    \vspace{-0.3cm}
    \caption{Approximation percent error and solve times for examples A and B. Solve times shown for objective $s$ are for one $\textup{FindApprox}(s,g_i)$ iteration.}
    \label{fig:timing_tiled}
    \vspace{-0.3cm}
\end{figure}

\vspace{-0.5cm}
\begin{equation*}
    \mathcal{X} = \{ x \in \mathbb{R}^2 \, | \, \begin{bmatrix} 1-16x_1x_2 & x_1 \\ x_1 & 1 - x_1^2 - x_2^2 \end{bmatrix} \succeq 0 \}.
\end{equation*}
% \begin{equation}
%     \mathcal{X} = \{ x \in \mathbb{R}^2 \, | \, 1-16x_1x_2 \geq 0, 1 - x_1^2 - x_2^2 \geq 0,  (1 - 16x_1x_2)(1 - x_1^2 - x_2^2) - x_1^2 \geq 0 \}
% \end{equation}
Using Algorithms \ref{alg:outer_Ker} and \ref{alg:inner_Ker} we find the kernel ($\KerOuter = \KerInner = \textup{conv}\{ \pm(-0.1752, 0.3335), \pm(0.1268, 0.2213) \}$) as shown in Figure \ref{fig:star_convex}. Figure \ref{fig:scale_example} (left) shows the 4th-order approximation obtained with Algorithm \ref{alg:inner_outer_X}. Figure \ref{fig:timing_tiled} shows the percent error as we increase the degree. Although each objective value (not shown) decreases monotonically with increasing degree, the percent error occasionally increases. This demonstrates the heuristic nature of the objectives for minimizing volume. 

\subsection{Discrete-time stabilizability region \textup{\cite{Henrion2012},\cite{Dabbene2017}}}
\begin{figure}
    \centering
    \vspace{0.15cm}
    %\vspace{-0.5cm}
    \includegraphics[width=0.44\textwidth]{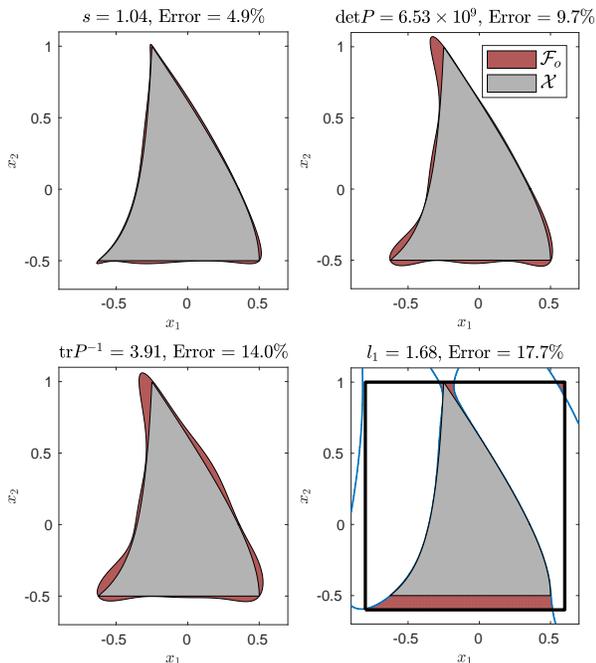}
    \vspace{-0.35cm}
    \caption{6th-order outer approximations of example B}
    \label{fig:Henrion2012_ex4p4_tiled}
    \vspace{-0.25cm}
\end{figure}

\vspace{-0.5cm}
\begin{align*}
    \mathcal{X} = \{ x \in \mathbb{R}^2 \, | \, 1 + 2x_2 \geq 0, 2 - 4x_1 - 3x_2 &\geq 0, \\ 10-28x_1-5x_2-24x_1x_2 - 18x_2^2 &\geq 0, \\
    1-x_2 - 8x_1^2 - 2x_1x_2 - x_2^2 - 8x_1^2x_2 - 6x_1x_2^2 &\geq 0    \}.
\end{align*}
The set contains the origin in its kernel. Figure \ref{fig:timing_tiled} shows the percent error for increasing degree. Figure \ref{fig:Henrion2012_ex4p4_tiled} shows the 6th-order approximations obtained with each objective. For the $l_1$ approximation we also show the bounding box from \cite{Dabbene2017}. 

% \subsection{General Non-Convex Set \textup{\cite{Cerone2012}}}
% \vspace{-0.5cm}
% \begin{equation*}
%     \mathcal{X} = \{ x \in \mathbb{R}^2 \, | \, (x_1-1)^2 + (x_2-1)^2 \leq 1, x_2 \leq 0.5x_1^2 \}
% \end{equation*}
% This set is star-convex but does not contain the origin in its kernel. In applying Algorithm \ref{alg:inner_outer_X} we translated the set to the Chebyshev center of its kernel $x^* = (1.39, 0.35)$.

\subsection{Convex Polytopes}
We generate 100 random convex polytopes in $\mathbb{R}^2$ with their Chebyshev center at the origin. We find outer approximations using the different objectives. Table \ref{tab:batch_cvx_2d} lists the number of times each objective obtained the smallest percent error relative to the other objectives for a given polytope. 

\begin{table}
%\vspace{1.75mm}
\caption{Instances In Which Objective Obtained Smallest Error}
\begin{center}
\begin{tabular}{c c c c c c}  
\hline 
  Deg. & \# Trials &  $s$ & $-\textup{det} P$ & $\textup{tr} P^{-1}$ & $l_1$\\
  \hline
  4 & 100 & 73 & 13 & 0 & 14 \\
  6 & 100 & 98 & 0 & 0 & 2
\end{tabular}

\label{tab:batch_cvx_2d}
\end{center}
\vspace{-1.75mm}
\end{table}

\textcolor{\editcolor}{
\subsection{Non-Star-Convex Set}
\vspace{-0.5cm}
\begin{align*}
    \mathcal{X} = \{ x \in \mathbb{R}^2 \, | \, r^2 \leq (x_1 - c)^2 + x_2^2 \leq 1, x_1 \leq c\}.
\end{align*}
Let $0 < r < c < 1$ so the origin is in the interior of the set. Figure \ref{fig:scale_example} shows the set for the case in which $c = 0.9$ and $r = 0.4$. Points $(c, \pm r) \in \partial \mathcal{X}$ yield cutting planes $x_2 \geq r$ and $x_2 \leq -r$ such that $\Ker \mathcal{X} = \emptyset$. Table \ref{tab:non_star_cvx} gives the outer approximation error for $c = 0.9$ and varying $r$.\footnote{The $l_1$ objective failed to improve upon the bounding box $\mathcal{B}$ supplied.}  For the scaling objective, we also report the objective value $s^\star$ and its lower bound $s_{lb}$.\footnote{The line segments connecting $(0,0)$ to $(c, \pm r)$ define the maximum lower bound on $s$ in Lemma \ref{lem:s_lowerbound}. It can be shown that $s_{lb} = \frac{\|p_2\|}{\|p_1\|}$ where $p_2 = (c,r), p_1 = (c + r\cos{\phi}, r\sin{\phi})$ and $\phi = \frac{\pi}{2} + 2\arctan{\frac{r}{c}}$.}
As $s_{lb}$ increases the percent error increases, confirming our heuristic is best suited to star-convex sets.
}

\begin{table}
\vspace{1.75mm}
\caption{Percent Error of Outer Approximations of Example E}
\begin{center}
\begin{tabular}{c c c c c}  
\hline
  r & Degree & $s (s^\star / s_{lb})$ & $-\textup{det} P$ & $\textup{tr} P^{-1}$\\ %[0.5ex] 
  \hline
  0.1 & 4 & 12.0 (1.096 / 1.025)    & 13.0  & 11.8 \\
  0.2 & 4 & 13.6 (1.104 / 1.104)   & 16.1  & 14.0 \\ 
  0.3 & 4 & 35.1 (1.250 / 1.250)    & 18.5  & 17.8 \\
  0.4 & 4 & 81.7 (1.492 / 1.492)   & 17.3  & 22.9 \\ 
 \hline
\end{tabular}
\label{tab:non_star_cvx}
\end{center}
\vspace{-0.75cm}
\end{table}
\textcolor{\editcolor}{
\subsection{Solver Performance}
 Figure \ref{fig:timing_tiled} shows the solve times for the various objectives on a logarithmic scale. Applied to a matrix $P \in \mathbb{R}^{m \times m}$, the $-\textup{det} P$ and $\textup{tr} P^{-1}$ objectives introduce a PSD matrix $H \in \mathbb{R}^{2m \times 2m}$ due to reformulations involving the exponential cone \cite{Mosek2017} and Schur complement \cite{Ahmadi2017} respectively. In contrast, the scaling $(s)$ and $l_1$ objectives work directly with $P$, yielding smaller semidefinite programs. The $l_1$ objective has the best computational performance. Due to the use of bisection, the total solve time for the scaling objective is an integer multiple of the time shown in Figure \ref{fig:timing_tiled}. Accounting for this, the scaling objective still remains competitive with the $-\textup{det} P$ and $\textup{tr} P^{-1}$ objectives.
 }
\subsection{Implementation Details}
YALMIP \cite{Lofberg2004} and MOSEK \cite{Mosek2017} were used to solve the SOS programs.\footnote{Supporting code will be released upon publication.} 
% In solving \eqref{opti:FindApprox}, we used polynomials $\lambda(x), \mu(x)$ with degree equal to that of the polynomial $f(x)$. 
\textcolor{\editcolor}{Volumes of non-star-convex sets were approximated by evaluating the indicator function over a discrete grid. Volumes of star-convex sets were approximated using numerical integration in polar coordinates.}

\section{CONCLUSIONS}
An algorithm for finding approximations of semialgebraic sets using sum-of-squares optimization was proposed. The algorithm relies on a novel objective which minimizes the scaling necessary to transform an inner approximation into an outer approximation of the set. Numerical examples demonstrated this objective often finds tighter approximations compared to existing heuristics when applied to star-convex sets. 
\textcolor{\editcolor}{Applied to non-star-convex sets, our proposed heuristic performs poorly. A promising direction to address this is through star-convex decompositions~\cite{Delanoue2006}. 
%Our approximation method may then be applied to the star-convex components. 
We leave this exploring this option for future work.}

\section*{ACKNOWLEDGEMENTS}
\textcolor{\editcolor}{
The author thanks Enrique Mallada and the anonymous reviewers for their valuable feedback.
}

\begin{appendix}
\subsection{Proof of Lemma 3}
% \begin{proof}
% Let $0 < a < b$ and $ax, bx \in \mathcal{X}$ for some $x \neq 0$. Let $b = 1$ without loss of generality. Let $tx \not\in \mathcal{X} \, \forall \, t \in (a, 1)$. 
% Let $F \subseteq \mathcal{X}$. Assume $s < \frac{1}{a}$ satisfies $\mathcal{X} \subseteq s\mathcal{F} \implies x \in s\mathcal{F} \implies \frac{x}{s} \in \mathcal{F}$. 
% Given $1 < s < \frac{1}{a} \implies a < \frac{1}{s} < 1 \implies \frac{x}{s} \not\in \mathcal{X} \implies \frac{x}{s} \not\in \mathcal{F}$ a contradiction. Thus $s \geq \frac{1}{a}$. The lower bound $s_{lb}$ follows by considering any $(a,x)$ satisfying the stated conditions.
% \end{proof}

% \textcolor{\editcolor}{Let $\mathcal{X}$ and $\mathcal{F}$ be compact sets in $\mathbb{R}^n$. Let $\mathcal{F} \subseteq \mathcal{X} \subseteq s^\star \mathcal{F}$ for some $s^\star > 1$. Let $0 \in \textup{int } \mathcal{F}$. 
% Let $x, sx \in \mathcal{X}$ and $tx \not\in \mathcal{X} \, \forall \, t \in (1,s)$ for some $s > 1, x \neq 0$. Then $s^\star \geq s$. }
\begin{proof}
\textcolor{\editcolor}{
Assume $1 < s^\star < s$ satisfies $\mathcal{F} \subseteq \mathcal{X} \subseteq s^\star \mathcal{F}$. Let $x, sx \in \mathcal{X}, x \neq 0$ such that $tx \not\in \mathcal{X} \, \forall \, t \in (1, s)$. Given $sx \in \mathcal{X} \implies sx \in s^\star\mathcal{F} \implies \frac{s}{s^\star}x \in \mathcal{F}$. 
However,  $1 < \frac{s}{s^\star} < s \implies \frac{s}{s^\star}x \not\in \mathcal{F}$, a contradiction.  Thus $s^\star \geq s$.
}
\end{proof}
\subsection{Proof of Lemma 4}
% If $\mathcal{F}$ is convex we can relate the scaling $s$ to the Hausdorff distance between the inner and outer approximation.
% \begin{lem}
% Let $\mathcal{F} \subseteq \mathbb{R}^n$ be a convex, compact set with $0 \in \textup{int} \mathcal{F}$. Let $s\mathcal{F}$ denote the scaled set with $s \geq 1$. 
% \begin{equation}
%     %d_H(s\mathcal{F},\mathcal{F}) = 
%     \underset{x \in s\mathcal{F}}{\textup {max}} \, \underset{y \in \mathcal{F}}{\textup {min}} \, |x-y| = (s-1)\cdot \underset{x \in \mathcal{F}}{\textup {max}} |x|
% \end{equation}
% \end{lem}
\begin{proof}
Recall the Hausdorff distance between two compact, convex sets can be written in terms of their support functions.
\begin{align}
%  \underset{x \in s\mathcal{F}}{\textup {max}} \, \underset{y \in \mathcal{F}}{\textup {min}} \, \|x-y\|_2 
d_H(s\mathcal{F},\mathcal{F})
 &=   \underset{c \in S^{n-1}}{\textup{max}} |\sigma_{s\mathcal{F}}(c) - \sigma_{\mathcal{F}}(c)| \\
&=\underset{c \in S^{n-1}}{\textup{max}} |s\sigma_{\mathcal{F}}(c) - \sigma_{\mathcal{F}}(c)| \\
&=(s-1)\cdot \underset{c \in S^{n-1}}{\textup{max}} \sigma_{\mathcal{F}}(c) \\
&=(s-1)\cdot \underset{x \in \mathcal{F}}{\textup {max}} \,\|x\|_2. 
\end{align}
\end{proof}

\subsection{Proof of Lemma 6}

\subsubsection{$\mathbf{\Ker ( \A \cap \B ) \supseteq \Ker\A \cap \Ker\B}$} 
Let $l(x,y) = \{ \lambda x + (1-\lambda)y \, | \, \lambda \in [0,1] \}$ for some $x \in \Ker\mathcal{A} \cap \Ker\mathcal{B} \textup{ and }  y \in \mathcal{A} \cap \mathcal{B}$. As $x \in \Ker\mathcal{A}, y \in \mathcal{A} \implies l(x,y) \subseteq \mathcal{A}$ and similarly, $x \in \Ker\mathcal{B}, y \in \mathcal{B} \implies l(x,y) \subseteq \mathcal{B}$, we see that $x \in \Ker(\mathcal{A} \cap \mathcal{B})$.\qedsymbol

\subsubsection{$\mathbf{\Ker ( \A \cup \B ) \supseteq \Ker\A \cap \Ker\B}$}
Let $l(x,y) = \{ \lambda x + (1-\lambda)y \, | \, \lambda \in [0,1] \}$ for some $x \in \Ker \mathcal{A} \cap \Ker \mathcal{B} \textup{ and }  y \in \mathcal{A} \cup \mathcal{B}$.  
For the case when $y \in \mathcal{A}$, then $x \in \Ker{A} \implies l(x,y) \subseteq \mathcal{A} \implies l(x,y) \subseteq \mathcal{A} \cup \mathcal{B}$. Similarly, for the case when $y \in \mathcal{B}$, then $x \in \Ker{B} \implies l(x,y) \in \mathcal{B} \implies l(x,y) \subseteq \mathcal{A} \cup \mathcal{B}$. Therefore $x \in \Ker (\mathcal{A} \cup \mathcal{B})$. \hfill \qedsymbol

\begin{remark}
\textcolor{\editcolor}{
Note that there is no relation between $\Ker (\mathcal{A} \cap \mathcal{B})$ and $\Ker ( \A \cup \B )$ in general. We gives examples in which one set is a subset of the other.\\
$\Ker (\mathcal{A} \cup \mathcal{B}) \supset \Ker (\mathcal{A} \cap \mathcal{B})$: Let $\mathcal{A} \setminus \mathcal{B} \neq \emptyset$ and $\mathcal{B} \setminus \mathcal{A} \neq \emptyset$. Let $\mathcal{A} \cup \mathcal{B}$ be a convex set. Then $\Ker(\mathcal{A} \cup \mathcal{B}) = \mathcal{A} \cup \mathcal{B} \supset (\mathcal{A} \cap \mathcal{B}) \supseteq \Ker(\mathcal{A} \cap \mathcal{B})$.} \\
\noindent \textcolor{\editcolor}{$\Ker (\mathcal{A} \cup \mathcal{B}) \subset \Ker (\mathcal{A} \cap \mathcal{B})$: Let $\mathcal{A}$ be a compact set that is not star-convex with non-empty interior. Let $\mathcal{B}$ be a non-empty convex set satisfying $\mathcal{B} \subset \mathcal{A}$. Then $\Ker(\mathcal{A} \cap \mathcal{B}) = \mathcal{B} \supset \emptyset = \Ker(\mathcal{A} \cup \mathcal{B})$.}
\end{remark}

\subsection{Proof of Lemma 5}
\newcommand{\ind}{{[n-1]}}
\newcommand{\indn}{{[n]}}
\newcommand{\kpt}{p} % Kernel point
\newcommand{\bpt}{q} % Boundary point
\nopagebreak

\begin{proof} 
\textcolor{\editcolor}{
For convenience, define the following:
\begin{equation*}
    \mathcal{H} :=  \{\kpt \,|\, \nabla g_i(\bpt)^T(\kpt - \bpt) \leq 0 \,\forall \, \bpt \in \partial \mathcal{X}_i, i \in [m] \}.
\end{equation*}
We show that $\Ker\mathcal{X} \subseteq \mathcal{H}$ and $\Ker\mathcal{X} \supseteq \mathcal{H}$ and therefore $\Ker\mathcal{X} = \mathcal{H}$.
}

\noindent $\Rightarrow (\Ker\mathcal{X} \subseteq \mathcal{H})$: 
Assume $\kpt \in \Ker\mathcal{X}$ but there exists a point $\bpt \in \partial\mathcal{X}_i$ for some $i \in [m]$ such that $\nabla g_i(\bpt)^T(\kpt - \bpt) > 0$.
Recall the definition of the directional derivative:
\begin{equation*}
\lim_{t\rightarrow0} \frac{g_i(t\kpt + (1-t)\bpt) - g_i(\bpt)}{t} =  \nabla g_i(\bpt)^T(\kpt - \bpt).    
\end{equation*}
Given $g_i(\bpt) = 1$ and $\nabla g_i(\bpt)^T(\kpt - \bpt) > 0$ implies there exists an open interval $t \in (0,\alpha), \alpha > 0$ in which $g_i(t\kpt + (1-t)\bpt) > 1$. The line segment over this open interval does not belong to $\mathcal{X}$. Thus $\kpt \not\in\Ker\mathcal{X}$, a contradiction.\newline
\textcolor{\editcolor}{
$\Leftarrow (\Ker\mathcal{X} \supseteq \mathcal{H})$: 
Let $\kpt \in \mathcal{H}$. Assume $\kpt \not\in \Ker\mathcal{X} \implies \exists \, \bpt \in \mathcal{X}$ such that $l(t) \not \in \mathcal{X}$ for some $t \in (0, 1]$ where $l(t) := t\kpt + (1-t)\bpt$.}
\footnote{\textcolor{\editcolor}{We have not yet shown that $\mathcal{H} \subseteq \mathcal{X}$ so we are not assuming $\kpt \in \mathcal{X}$.}}
\textcolor{\editcolor}{
As $\mathcal{X}$ is compact, $l(t) \not \in \mathcal{X} \implies g_i(l(t)) > 1$ for some $i \in [m]$ and open interval $t \in (a,b)$ satisfying $0 \leq a < b$ with $a < 1$.
Without loss of generality, let $a = 0$ such that $\bpt \in \partial \mathcal{X}_i$ and $g_i(l(0)) = 1$. Applying the definition of the directional derivative yields:
\begin{equation*}
\lim_{t\rightarrow0} \frac{g_i(l(t)) - g_i(l(0))}{t} =  \nabla g_i(\bpt)^T(\kpt - \bpt).    
\end{equation*}
The left-hand side of this relation is non-negative. The right-hand side is non-positive per the definition of $\mathcal{H}$. Thus both sides must equal zero. As $\nabla g_i(\bpt) \neq 0$, this implies 
\begin{equation}\label{eqn:perp}
(\kpt - \bpt) \perp \nabla g_i(\bpt).
\end{equation}
}
\textcolor{\editcolor}{
Assume w.l.o.g. that $\nabla g_i(\bpt)$ is aligned with coordinate $n$:
\begin{equation}\label{eqn:y_align}
    \nabla g_i(\bpt) = \begin{bmatrix} 0_{n-1}^T && r \end{bmatrix}^T, r > 0.
\end{equation}
}
\textcolor{\editcolor}{
If this does not hold we can introduce an appropriate change of variables. Together, \eqref{eqn:perp} and \eqref{eqn:y_align} $\implies (\kpt_n - \bpt_n)r = 0 \implies l_n(t) = \bpt_n$. From this we have
\begin{equation}
    l(t) = \begin{bmatrix}
    t\kpt_\ind + (1-t)\bpt_\ind \\
    \bpt_n
    \end{bmatrix}.
\end{equation}
}
\textcolor{\editcolor}{
Define the following parameterized curve $\phi: \mathbb{R} \rightarrow \mathbb{R}^n$ which moves along the boundary $g_i(x) = 1$, starting from $\kpt$:
\begin{equation}
    \phi(t) = \begin{bmatrix}
    t\kpt_\ind + (1-t)\bpt_\ind \\
    h(t\kpt_\ind + (1-t)\bpt_\ind)
    \end{bmatrix}.
\end{equation}
}
\textcolor{\editcolor}{
Given $\frac{\partial g_i}{\partial x_n}(\bpt) \neq 0$, from the implicit function theorem there exists an open set $U \subset \mathbb{R}^{n-1}$ with $\bpt_\ind \in U$ and $C^1$ function $h: U \rightarrow \mathbb{R}$ such that $h(\bpt_\ind) = \bpt_n$ and $g_i(x_\ind,h(x_\ind)) = 1$ for all $x_\ind \in U$. Here we are restricting coordinates $x_\ind$ to the line segment parameterized by $t$. Thus $g_i(\phi(t)) = 1$ for all $t$ such that $\phi_\ind(t) \in U$. Let $t \in (-c,\,d), c > 0, d > 0$ denote this interval.
}
\textcolor{\editcolor}{
The line $l(t)$ and curve $\phi(t)$ only differ in coordinate $n$. Given $g_i(l(t)) > 1, t \in (0,b)$ and $g_i(\phi(t)) = 1, t \in (-c, d) \implies \bpt_n \neq \phi_n(t) \, \forall \, t \in (0, \textup{min}(b,d))$. }
\textcolor{\editcolor}{
Given $\frac{\partial{g_i}}{\partial x_n}(\bpt) > 0 \implies \frac{\partial{g_i}}{\partial x_n} > 0$ for some open ball around $\bpt$ as $g_i$ is smooth.  Assuming $\phi_n(t) > \bpt_n \implies g_i(\phi(t)) > g_i(l(t)) > 1$ for points sufficiently close to $\bpt$, a contradiction. Thus $\phi_n(t) < \bpt_n$ for some interval $t \in (0, e), e > 0$. From this we have
}
\textcolor{\editcolor}{
\begin{equation}\label{eqn:dy_neg2}
    \frac{\partial g_i(\phi(t))}{\partial x_n}(\bpt_n - \phi_n(t)) > 0, \forall t \in (0, e).
\end{equation}
}
\textcolor{\editcolor}{
Given $\bpt_n = \phi_n(0) > \phi_n(t)$ for some interval $t \in (0, e)$, by the mean value theorem there exists $t_\star \in (0,e)$ such that $\frac{d\phi_n}{dt}(t_\star) < 0$. This yields the following relation:
\begin{equation}\label{eqn:dy_neg}
    \frac{\partial g_i(\phi(t_\star))}{\partial x_n} \frac{d\phi_{n}(t_\star)}{dt} < 0. 
\end{equation}
}
\textcolor{\editcolor}{
Given $g_i(\phi(t)) = 1 \, \forall \, t \in (-c, d) \implies \frac{dg_i}{dt}(\phi(t)) = 0$. We expand this at the point $t_\star$ obtaining
\begin{equation}\label{eqn:dg_zero}
\begin{split}
    0 &= \frac{\partial g_i(\phi(t_\star))}{\partial x_\ind}^T  \frac{d\phi_\ind(t_\star)}{dt} + 
    \frac{\partial g_i(\phi(t_\star))}{\partial x_n} \frac{d\phi_{n}(t_\star)}{dt} \\
    &= \frac{\partial g_i(\phi(t_\star))}{\partial x_\ind}^T(\kpt_\ind - \bpt_\ind) + 
    \frac{\partial g_i(\phi(t_\star))}{\partial x_n} \frac{d\phi_{n}(t_\star)}{dt}.
    %\nabla_w g_i(\phi(t_\star))^T(w_k - w_b) + \nabla_y g_i(\phi(t_\star))\frac{dh}{dt}(t_\star) = 0
\end{split}
\end{equation}
From equations \eqref{eqn:dy_neg} and \eqref{eqn:dg_zero} we obtain
\begin{equation}\label{eqn:dw_pos}
    \frac{\partial g_i(\phi(t_\star))}{\partial x_\ind}^T(\kpt_\ind - \bpt_\ind) > 0.
\end{equation}
Finally, we evaluate the stated constraint on $\kpt \in \mathcal{H}$ at the boundary point $\phi(t_\star)$ giving 
\begin{equation}
    \begin{split}
        \nabla g_i(\phi(t_\star))^T(\kpt - \phi(t_\star))  = 
        \frac{\partial g_i(\phi(t_\star))}{\partial x_n} (\kpt_n - \phi_n (t_\star)) \, + \, \\
        \frac{\partial g_i(\phi(t_\star))}{\partial x_\ind}^T(\kpt_\ind - \bpt_\ind)(1-t_\star). \\        
    \end{split}
\end{equation}
From \eqref{eqn:dy_neg2} and \eqref{eqn:dg_zero} and noting that $(1-t_\star) > 0$ and $q_n = p_n$ gives
\begin{equation}
    \nabla g_i(\phi(t_\star))^T(\kpt - \phi(t_\star))  > 0.
\end{equation}
Thus $\kpt \not \in \mathcal{H}$, a contradiction.
}
\end{proof}

\end{appendix}
\addtolength{\textheight}{-12cm}   % This command serves to balance the column lengths

%\nocite{*}

\bibliography{references}
\bibliographystyle{ieeetr}

\end{document}